\documentclass[12pt]{amsart}
\usepackage{amssymb,amscd,verbatim}

\swapnumbers

\newcommand{\Hom}{\operatorname{Hom}}      
\newcommand{\Ext}{\operatorname{Ext}}      
\newcommand{\DM}{{\rm DM}}          
\newcommand{\M}{\mathcal{M}_1}   

\newcommand{\car}{\operatorname{char}}

\newcommand{\ihom}{{\rm\underline{Hom}}}  

\renewcommand{\P}{\mathbb{P}}   
\renewcommand{\L}{\mathbb{L}}  
\newcommand{\Aff}{\mathbb{A}}   

\newcommand{\C}{\mathbb{C}}     
\newcommand{\Q}{\mathbb{Q}}     
\newcommand{\Z}{\mathbb{Z}}     

\newcommand{\G}{\mathbb{G}}     
\newcommand{\HH}{\mathbb{H}}    
\newcommand{\EExt}{{\rm \mathbb{E}xt}} 

\newcommand{\im}{\operatorname{Im}}        
\renewcommand{\ker}{\operatorname{Ker}}  
\newcommand{\gr}{\operatorname{gr}}        
\newcommand{\Pic}{{\rm Pic}}     
\newcommand{\Alb}{{\rm Alb}}     
\newcommand{\RPic}{{\rm RPic}}     
\newcommand{\LAlb}{{\rm LAlb}}     
\newcommand{\LA}[1]{\mbox{${\rm L}_{#1}{\rm Alb}$}}
\newcommand{\RA}[1]{\mbox{${\rm R}^{#1}{\rm Pic}$}}

\newcommand{\Tot}{{\rm Tot}}     

\newcommand{\by}[1]{\stackrel{#1}{\rightarrow}}
\newcommand{\longby}[1]{\stackrel{#1}{\longrightarrow}}

\newcommand{\df}{\mbox{\,${:=}$}\,}
\newcommand{\ie}{{\it i.e.\/},\ }
\newcommand{\cf}{{\it cf.\/}\ }
\newcommand{\eg}{{\it e.g.\/},\ }
\newcommand{\et} {\mbox{\scriptsize{\rm {\'e}t}}}

\newcommand{\fr}{\mbox{\scriptsize{\rm fr}}}

\newcommand{\eff}{\mbox{\scriptsize{\rm eff}}}
\newcommand{\gm}{\mbox{\scriptsize{\rm gm}}}

\renewcommand{\bar}{\overline}
\newcommand{\into}{\hookrightarrow}
\renewcommand{\implies}{\mbox{$\Rightarrow$}}
\newcommand{\sZ}{\mbox{\scriptsize{$\Z$}}}   
\newcommand{\sC}{\mbox{\scriptsize{$\C$}}}   
\newcommand{\sQ}{\mbox{\scriptsize{$\Q$}}}   

\newcommand{\liminv}[1]{\mathop{\rm
lim}_{\buildrel\longleftarrow\over{#1}}}
\newcommand{\onto}{\mbox{$\to\!\!\!\!\to$}}

\newcommand{\boxtensor}{\def\boxtimesten{\Box\kern-7.59pt\raise1.2pt
\hbox{$\times$} }}                                  

\newcounter{elno}                   

\newcommand{\cC}{\mathcal{C}}
\newcommand{\cD}{\mathcal{D}}

\newcommand{\cF}{\mathcal{F}}

\newcommand{\cM}{\mathcal{M}}

\newcommand{\cO}{\mathcal{O}}

\newcommand{\cV}{\mathcal{V}}

\newcommand{\cZ}{\mathcal{Z}}

\begin{document}

\title{A pamphlet on Motivic Cohomology}
\author{Luca Barbieri-Viale}
\address{Dipartimento di Matematica Pura e Applicata, Universit\`a degli Studi di Padova\\ Via G. Belzoni, 7\\Padova -- I-35131\\ Italy}
\email{barbieri@math.unipd.it}
\date{July 15, 2005}
\begin{abstract} 
This is an ``elementary'' introduction to the conjectural theory of {\it motives}\, along the lines indicated by Grothendieck. We further quote recent developments, also presenting some advances due to Voevodsky, and applications to the study of algebraic cycles and differential forms.
\end{abstract}
\maketitle
\hfill {\small \it Ad astra per aspera}\\

In order to fix a general ``philosophical'' framework of what  a theory of {\em motives}\, looks like let's draw the following picture\\
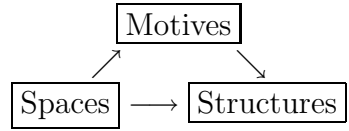
$$\begin{array}{c} 
\fbox{{\rm Motives}}\\
\displaystyle
\nearrow\hspace*{1.5cm} \searrow\\
\fbox{{\rm Spaces}} \longby{} \fbox{{\rm Structures}}\\
\end{array}$$
with the following explanation:
\begin{itemize}
\item {\bf Spaces $ \implies $ Structures}: we usually go from spaces to structures associating several kind of {\em invariants} of the {\em shape} we are investigating\footnote{The most elementary is the dimension but we here imagine quite many other structures corresponding to homotopical, topological, differential or algebraic invariants which are linked to several visions of our spaces embodied in topological, differential and algebraic varieties or manifolds.} and call them {\em cohomology} theories of our spaces along with the structures they carry on;
\item {\bf Spaces $ \implies $ Motives}: when a concept of space is fixed and a concept of cohomology theory is involved we then seek for a {\em universal} cohomology theory along with its structure of {\it motive} and call it {\it motivic cohomology};
\item  {\bf Motives $ \implies $ Structures}: such a {\it motive} of the shape will be the finest invariant or structure associated to a space and would have several {\em realizations} yielding the various cohomology theories.\footnote{Exactly as in the common sense we here just think of a musical or visual motive being further realized with several different instruments.}\\ \end{itemize}

The required universality will be rendering {\it all} cohomology theories as realizations of {\it one} motivic cohomology. If we then prove a theorem for {\it motives} we prove several theorems at once, corresponding to the various realizations. This is the main motivation for looking at motivic cohomology in algebraic geometry and its existence will answers to the question or ``mistery'' of parenthood between {\em a priori} different cohomology theories\footnote{For example, we here have in mind,  in algebraic geometry, the so called $\ell$-adic cohomologies as the prime $\ell$ varies, provided by Grothendieck \'etale cohomology, see \ref{ex} below.} having deep similarities, providing related informations and interweaving arithmetic and geometry. 

For example, regarding the picture above in topology where ``spaces'' are CW-complexes and the cohomology theories are those satisfying wedge and Mayer-Vietoris axioms we obtain the so called spectra as ``motives'' and the Brown representability theorem is granting the universality. After the revolution due to Poincar\'e's {\it analysis situs} and the spread out of algebraic topology in geometry, Grothendieck's concept of {\it motive} is a glimpse of the future of algebraic geometry.

In order to keep the following presentation short and of interest to non  specialists we have to be quite rough and sometimes vague (we apologize to the expert reader).

I would like to thank Y. Andr\'e and M. Saito for some useful comments on the preliminary version of this note.

\section{Pure motives}
\subsection{Cohomology theory} Let $k$ be an algebraically closed field, \eg $k=\C$ the complex numbers. Usually, a cohomology theory on the category of non-singular projective varieties $\cV_k$ is   
 $$X\mapsto H^*(X)$$
associating to each $X\in \cV_k$ a graded group $H^*(X)\df \sum_i H^i(X)$ with a {\it cup product} $$H^i(X)\times H^j(X)\to H^{i+j}(X)$$
denoted $(\alpha, \beta)\leadsto \alpha \cdot \beta$ and 
such that every morphism $f: X\to Y$ yields $f^*: H^*(Y)\to H^*(X)$
compatibly with the compositions and identities, \ie yielding a functor.
Moreover, $H^*(X)$ is equipped with a suitable structure.

\subsubsection{ } \label{cl} A key requirement, linking geometry to topology, is the following {\it cycle map}. If $X \in \cV_k$ then there is a group homomorphism
$$c\ell^i_X : \cZ^i (X) \to H^{2i}(X)$$
where $\cZ^i (X) \df$ ``algebraic cycles'' of codimension $i$.
The algebraic cycles are finite ``linear combinations'' of closed sub-varieties, \ie $\cZ^i (X)$ is the free abelian group generated by closed sub-varieties of codimension $i$ on $X$.

\subsubsection{ } \label{ax} When $H^i(X)$ are actually finite dimensional vector spaces for some coefficients field $K$, \eg $K=\Q$ the rational numbers, the wonderful properties\footnote{These properties are peculiar of the so called Weil cohomology theories.} we may have are the following. For $X\in \cV_k$ set the algebraic dimension $n\df \dim (X)$ (\eg here a ``curve'' is a Riemann surface) and 
\begin{enumerate} 
\item {\it Poincar\'e duality}: we have an ``orientation'' such that $H^{2n}(X)\cong K$ and there is a ``duality'' induced by the cup product
$$H^i(X)^{\vee}\cong H^{2n-i}(X)$$
where $H^i(X)^{\vee}\df \Hom (H^i(X),K)$ is the dual $K$-vector space;\\

\item {\it K\"unneth formula}: if $X, Y \in \cV_k$ then
$$H^k(X\times Y)\cong \sum_{i+j=k} H^i(X)\otimes H^j(Y)$$

\item {\it Hard Lefschetz}: for $H$ a smooth hyperplane section of $X$ define the Lefschetz operator $L: H^i(X) \to H^{i+2}(X)$
 by $L (\alpha) \df \alpha \cdot c\ell (H)$. For $i\leq n $ the iterated Lefschetz operator induces
 $$H^i(X)\cong H^{2n-i}(X)$$
\end{enumerate}

Let $A^* (X)$ denote the sub-group of $H^{2*}(X)$ given by algebraic cycles, \ie by the image of $c\ell^*$ above. This $A^* (X)$ is the so called ``algebraic part'' of the cohomology and is conjecturally independent of a particular choice of a cohomology theory $H^* (X)$ satisfying the properties listed above.

\subsubsection{ } \label{ex} The following examples are the main examples of cohomology theories which do have the properties listed in \ref{cl} and \ref{ax} 
\begin{enumerate}
\item {\it Singular cohomology}: for $k=\C$ an algebraic variety $X\in \cV_k$ can be regarded as a complex manifold $X_{\rm an}$ and we can set
$$H^*(X)\df H^*(X_{\rm an}, \Z)\otimes K$$
taking the classical $K$-linearized singular cohomology of the associated analytic space;
\item {\it \'Etale cohomology}: for $k=\bar k$ of arbitrary characteristics, \eg the algebraic closure of a finite field, an algebraic variety $X\in \cV_k$ can be regarded to carry on the \'etale topology\footnote{Grothendieck topologies are allowing certain ``opens'' which are not sub-sets. One abstractly defines some families of maps to be  ``coverings''  in a category obtaining a so called ``site'' and the corresponding Grothendieck ``topos'' of sheaves. The \'etale site is finer than the classical Zariski site and was introduced by Grothendieck in order to detect arithmetic properties not detected by the Zariski topology.} denoted by $X_{\et}$
and we can set
$$H^*(X)\df \liminv{\nu}H^*(X_{\et}, \Z/\ell^{\nu})\otimes K$$
where usually $K=\Q_{\ell}$ (these are $\Z_{\ell}$-modules and $\ell$ is a prime different from the characteristic of $k$).
When $k=\C$ the \'etale topology is related to the usual analytic topology by a ``continuous'' map $X_{\rm an}\to X_{\et} $ yielding a ``comparison'' isomorphism
$H^*(X_{\et},\Z/{\ell}^{\nu})\cong H^*(X_{\rm an},\Z/{\ell}^{\nu})$.
\end{enumerate}

\subsection{Grothendieck motives} From \ref{ax} we have
$$H^*(X\times Y)= H^*(X)\otimes H^*(Y)  = H^*(X)^{\vee}\otimes H^*(Y) = \Hom (H^*(X), H^*(Y) ) $$
The general principle suggested by this identification is that {\it any} linear operator $\Psi: H^*(X)\to H^*(Y) $ on the cohomology, which is of algebraic nature, will be possibly defined by an algebraic cycle $\psi\in A^*(X\times Y)\subseteq H^*(X\times Y)$ that is independent of the cohomology theory. It would be the case if $A^*(X)$ (as a $\Q$-vector space) itself provides such a cohomology theory and this is roughly the content of the Standard Conjectures \cite{SG}:
\begin{itemize}
\item the {\it Lefschetz Standard Conjecture} claims that the Lefschetz operator $L$ induces an isomorphism on $A^*(X)$;
\item  the {\it Hodge Standard Conjecture} claims that the cup product on ``primitive'' classes in $A^*(X)$ is positive definite.\footnote{Primitive classes are those in the kernel of the Lefschetz operator. In characteristic zero, this conjecture is true for \'etale cohomology.}
\end{itemize}
Moreover, $ A^*(X\times Y)$ will then provide homomorphisms in the category of motives as follows. 
\subsubsection{ } \label{cyc} Observe that $A^*(X)\df (\cZ^*(X)/\equiv_{\rm hom})\otimes \Q$ where $\equiv_{\rm hom}$ is the ``cohomological'' equivalence relation given by the kernel of the cycle map \ref{cl}. Actually, there are several ``adequate'' relations between algebraic cycles, \eg two sub-varieties $Z_0$ and $Z_1$ of $X$ are {\it rationally} equivalent if they appear in a family $\{Z_t\}$ parametrized by $\P^1$. This  rational equivalence $\equiv_{\rm rat}$
is  ``minimal'' such that the intersection becomes a product, \ie if we set  the so called Chow groups $CH^*(X)\df \cZ^*(X)/\equiv_{\rm rat}$ we have $c\ell^i_X : CH^i(X)\to H^{2i}(X)$ and
$$CH^i(X)\times CH^j(X)\to CH^{i+j}(X)$$
providing the ``intersection'' product\footnote{However, $CH^i(X)\otimes K$ are not finite dimensional!} compatibly with \ref{cl}. Assuming the Standard Conjectures the resulting ($\Q$-linearized) pairing $$A^i(X)\times A^{n-i}(X)\to \Q$$
is non degenerate\footnote{In characteristic zero, this property is equivalent to the named conjectures.}  and therefore $\equiv_{\rm hom}$ is simply given by the so called {\it numerical} equivalence $\equiv_{\rm num}$ provided by the intersection pairing (independently  of $H^*(X)$). 

\subsubsection{ } \label{cor} A cycle in the group $\cZ^*(X\times Y)$ is called a {\it correspondence} from $X$ to $Y$ (with rational coefficients when $\Q$-linearized).
For any ``adequate'' relation $\equiv$ we obtain a category $\cC\cV_k$ of correspondences given by $X, Y\in\cV_k$ and $\psi: X\leadsto Y$ where $\psi\in \cZ^n(X\times Y)/\equiv\df \Hom_{\cC\cV_k}(X, Y)$. The composition of  $\psi: X\leadsto Y$ and  $\phi: Y\leadsto Z$  is simply given by pulling back $\psi$ and $\phi$  in $X\times Y\times Z$ and pushing forward  on $X\times Z$ their intersection product.

\subsubsection{ } \label{mot} The category $\cM^{\eff}_{\equiv}$ of {\it effective motives} is then formally obtained as follows. The objects are pairs $(X, p)$ with $X\in\cV_k$ and $p:X\leadsto X$ is a {\it projector} so that $p^2=p$. A morphism from $M = (X, p)$ to $N = (Y, q)$ is a 
correspondence $X\leadsto Y$ which is compatible with $p$ and $q$. We then obtain
the motive of an algebraic variety
$$X \mapsto M (X)\df (X,\Delta)$$
associating $X\in\cV_k$ with the diagonal cycle $\Delta : X\to X\times X$ which is the identity $X\leadsto X$. We then associate a map $f:X\to Y$ to its graph\footnote{It is usually taken the transposed of its graph in order to have a contravariant theory.} $\Gamma_f \subseteq X\times Y$ in such a way that $M: \cV_k\to \cM^{\eff}_{\equiv}$. 

The category $ \cM^{\eff}_{\equiv}$ has direct sums $\oplus$ and tensor products $\otimes$ given by disjoint unions
and products of varieties respectively. Furthermore, every projector $p$ on $M$ has a kernel and provides a decomposition $M = \ker (p) \oplus \ker (1-p)$.

\subsubsection{ } \label{lef} Let $e\in X$ be a point and consider the projectors  $X\leadsto X$ defined by $p_0\df e\times X$ and $p_{2n}\df X\times e$. It is easy  to see that $\Delta - p_{0}- p_{2n}$ is still a projector. For a curve $X$ we set  $M^0 (X)\df (X,p_0)$, $M^2(X)\df (X, p_{2})$ and $M^1(X) \df (X, p_1)$ where $p_1\df \Delta - p_{0}- p_{2}$ providing a decomposition
$$M (X) = M^0 (X) \oplus M^1  (X) \oplus  M^2 (X)$$
In particular, we set $\L \df M^2 (\P^1)$, the {\it Lefschetz} motive, and the category $\cM_k$ of Grothendieck motives\footnote{Here the equivalences ${\equiv}_{\rm hom}$ and   ${\equiv}_{\rm num}$ in \ref{cyc} are considered with rational coefficients and are assumed to coincide so that we obtain a unique $\Q$-linear abelian category.} over $k$ is obtained by formally adding  to  $\cM_{{\equiv}_{\rm hom}}^{\eff}$ the tensor product inverse $\L^{-1}$ called the {\it Tate} motive.\footnote{This is a technical device in order to make motivic Tate ``twists'' in Poincar\'e duality. For the sake of simplicity, we have not mentioned Tate twists in \ref{ax}.(1) but they will  appear in \ref{mhs} below. Tate twists are usually omitted as far as they yield an automorphism.}
\subsubsection{ } \label{gal} Grothendieck Standard Conjectures \cite{SG} then also grant that  the projector provided by the composition of the projection and the inclusion 
$$H^*(X)\onto H^i(X)\into H^*(X)$$
is induced by an algebraic cycle $p_i$ and, therefore, the above decomposition holds true (with $\Q$-coefficients)  in higher dimension in such a way that $M (X) = M^*(X)= \sum_i M^i(X)$ where $M^i(X) \df (X, p_i)$ is\footnote{Since these projectors $p_i$ are also the K\"unneth components of the diagonal cycle {\it via} \ref{ax}.(2), Murre \cite{MU} conjectures that such a decomposition can be lifted modulo rational equivalence and the resulting Chow-K\"unneth decomposition will provide a filtration on the Chow groups, see \ref{filt}.} formally a cohomology theory, the named {\sf motivic cohomology} of smooth projective varieties.\footnote{Note that the resulting objects are not abelian groups as in the usual definition \ref{cl} of a cohomology theory.} We then have that {\it any} such a cohomology theory $H^*$ (as described in \ref{ax}) factors
$$\begin{array}{c} 
\cM_{k}\\
\displaystyle
 {}_{ M^*}\nearrow \hspace*{0.4cm} \searrow\hspace*{0.5cm}\\
\cV_k\longby{H^*} \mbox{ ?}\\
\end{array}$$
This $\Q$-linear category $\cM_{k}$ would be an abelian\footnote{An abelian category is an  ``abstract'' version of the category of abelian groups in such a way that we can deal with exact sequences.} semi-simple\footnote{Semisimple means that every object has finite length and any sub-object a complement, therefore, all exact sequences split exactly like for finite dimensional $K$-vector spaces.} tensor category which is Tannakian which means that $\cM_{k}$ will be also equivalent to a category of representations for a pro-algebraic group: the motivic Galois group. We strongly recommend \cite[Part I]{A} for an extended  valuable orientation and for references on these matters. We also refer to \cite{M} and \cite{MU} for  further details.

\section{Mixed motives}
\subsection{Mixed Hodge structures} For $X$ a complex algebraic manifold we have De Rham cohomology $H^*_{DR}(X) \cong H^*(X_{\rm an}, \Z)\otimes \C$ embodied in singular cohomology (see \ref{ex}.(1)) and the De Rham decomposition
$$H^r_{DR}(X) = \sum_{p+q=r}H^p(X, \Omega^q_X)$$
providing the Hodge filtration: singular cohomology is supporting a bigraded structure called {\it Hodge structure}.  The singular cohomology groups $H^r(X_{\rm an}, \Z)$  of {\it any} complex algebraic variety, which can be open and singular, are endowed with a more sophisticated structure: a {\it mixed} Hodge structure, discovered by Deligne \cite{D}.

\subsubsection{ } \label{mhs} It is abstractly defined as a triple $H \df  (H_{\sZ}, W, F)$ where $H_{\sZ}$ is a finitely generated abelian group, \eg $H^r(X_{\rm an}, \Z)$,  the so called {\it weight} filtration $W_j$ is a finite increasing filtration on $H_{\sZ}\otimes\Q$ and the {\it Hodge} filtration $F^i$ is a finite decreasing filtration on $H_{\sZ}\otimes\C$ such that $W, F$ and its conjugate $\bar F$ is a system of ``opposed'' filtrations: there is a canonical decomposition  $$\gr_r^{W}(H_{\sQ}) \otimes \C =
\sum_{p+q=r}^{}
H^{p,q}$$ where $H^{p,q} \df F^p\cap \bar{F}^q$ and conversely.
When $\gr_r^{W}\neq 0$ for just a single $r$ we say that $H$ is {\it pure} of weight $r$, \eg this is the case of $H^r(X_{\rm an}, \Z)$ when $X$ is smooth and projective. Define $\Z (1)$ as the mixed Hodge structure on $2\pi\sqrt{-1} \Z$ pure of weight $-2$ purely of Hodge type $(-1,-1)$. The Tate twist $\Z (n)\df \Z (1)^{\otimes n}$ is then pure 
of weight $-2n$ and purely of Hodge type $(-n,-n)$.
We may also define the {\it level}\,  by $\ell (H)\df \max \{|p-q|:H^{p,q}\neq 0\}$.
The abelian category of $\Z$-mixed Hodge structures ${\rm MHS}$ has objects $H$ as above and morphisms $\varphi : H \to H'$ preserving the filtrations. The
kernel (resp. the cokernel) of a morphism $\varphi : H \to H'$ has
underlying  $\Q$  and $\C$-vector spaces the kernels (resp. the cokernels)
of $\varphi_{\sQ}$ and $\varphi_{\sC}$ with induced filtrations and any
morphism is strictly compatible with the filtrations.\footnote{The functors
$\gr_{W}$ and  $\gr_{F}$ are exacts.}

\subsubsection{ } \label{ext}
Let $Z\into X$ be a closed sub-variety  of the complex
algebraic variety $X$.  The relative singular cohomology $H^*_Z(X)\df H^*(X, X-Z;\Z)$ (\cf \cite[8.2.2 \&  8.3.8]{D}) carry on a mixed
Hodge structure fitting into a long exact sequence in ${\rm MHS}$
$$\cdots \to H_Z^i(X)\to H^i(X)\to
H^i(X-Z)\to H_Z^{i+1}(X)\to \cdots $$
We then get a refined cohomology theory (extending \ref{cl} and \ref{ax})
$$Z\subseteq X \mapsto H_Z^*(X)$$
with values in the abelian category of mixed Hodge structures.\footnote{It is a so called Poincar\'e duality theory with supports which is appropriate for algebraic cycles, \ie we also have a Borel-Moore homology theory $H_*^{BM}(X)$ and there is a ``cap'' product $H_Z^i(X)\otimes H_j^{BM}(X)\to H_{j-i}^{BM}(Z)$ such that,  when $X$ is smooth, $H_Z^i(X)\cong H_{2n-i}^{BM}(Z)$ by capping with the ``fundamental class''  which is a generator of $H_{2n}^{BM}(X)$ for $X$ smooth $n$-dimensional.}
For example, let $X-Z$ be a smooth curve obtained by removing a finite set  $Z$ of closed points from $X$ a smooth compact curve: we then get the mixed $H^1(X-Z)$ as an ``extension'' by pure objects  as above. In general, in the category ${\rm MHS}$ there are non trivial extensions and the $\Ext$ in {\rm MHS} naturally provides deep geometrical informations.

\subsubsection{ } \label{hdg} Let  $Z\subseteq X$ be a closed sub-variety of codimension $i$ in a projective algebraic manifold $X$. It is not difficult to see that the cycle class $c\ell (Z)$  in De Rham cohomology belongs to $H^i(X, \Omega^i_X)\subseteq H^{2i}_{DR}(X)$. 
Those singular cohomology  classes in $H^{2i}(X)$ landing in  $H^i(X, \Omega^i_X)$ are called {\it Hodge cycles}. Let $H^{i,i}_{\sQ}(X)\subseteq H^{2i}(X_{\rm an}, \Q)$ denote the rational Hodge cycles. We then have $A^i(X)\subseteq H^{i,i}_{\sQ}(X)$ and 
\begin{itemize}
\item the {\em  Hodge conjecture}  claims that the equality $A^i(X) = H^{i,i}_{\sQ}(X)$ holds.
\end{itemize}

Since the Lefschetz operator induces $ H^{i,i}_{\sQ}(X)\cong  H^{n-i,n-i}_{\sQ}(X)$ the Hodge conjecture implies the Lefschetz Standard conjecture.\footnote{The Standard Conjectures in characteristic zero are weaker than the Hodge Conjecture.}

The main historical evidence for this quite classical conjectures comes from the Lefschetz theorem for $(1,1)$-classes, \ie from the complex exponential $\exp : \C \to \C^*$ yielding an exact sequence
$$ CH^1 (X)\by{c\ell} H^2(X_{\rm an},\Z)\to H^2(X,\cO_X)$$
which characterize the image of $c\ell$ since $H^2_{DR}(X)/F^1 = H^2(X,\cO_X)$ and 
$H^{1,1}_{\sZ}(X) = H^{2}(X_{\rm an}, \Z)\cap F^1H^{2}_{DR}(X)$. Note that here, since $X$ is smooth, Weil divisors coincide with Cartier divisors so that $ CH^1 (X)= \Pic (X)$ and $\Pic (X) = \Pic (X_{\rm an})$ since $X$ is proper, granting also that $H^*(X, \cO_X)= H^*(X_{\rm an},\cO_{X_{\rm an}})$.

\subsubsection{ } \label{ghdg} In general, the Grothendieck coniveau or arithmetic filtration (\cf \cite{GH}) $$N^iH^j(X) \df \bigcup_{{\rm codim}_XZ\geq i} \ker (H^j(X) \to  H^j(X-Z))$$ yields a filtration by Hodge sub-structures of $H^j(X)$. We have that
$$N^iH^j(X, \Q)\subseteq H^j(X, \Q)\cap F^iH^j_{DR}(X)$$
where $N^iH^j(X, \Q )$ is of level $\leq j-2i$ and we have \cite{GH}:
\begin{itemize}
\item  the  {\em Grothendieck-Hodge conjecture} claims that the left hand side of the inclusion  above is the largest sub-structure of level $\leq j-2i$ of the right hand side.
\end{itemize}
Note that for $Z\subset X$, $Z\in  \cZ^i(X)$, then $c\ell (Z)\in H^{2i}(X)$ is exactly the image of $H^{2i}_Z(X)$,  it vanishes in $H^{2i}(X-Z)$ and $H^{i,i}_{\sQ}(X) = H^{2i}(X_{\rm an}, \Q)\cap F^iH^{2i}_{DR}(X)$. In general, for $Z$ smooth $H_Z^j(X) =H^{j-2i}(Z)$ and the exact sequence
$$H^{j-2i}(Z)\to H^i(X)\to H^i(X-Z)$$
is clarifying a bit more the assertion above.\footnote{We can find counterexamples to naive formulations of these assertions for singular varieties but there are perfectly suitable reformulations, see \cite{BV} for the precise statements, making use of the theory of 1-motives, see \ref{1-mot} below.}

\subsection{Voevodsky motives} By dealing with compact algebraic manifolds Grothendieck Standard Conjectures grant ``universal'' properties of existing categories of motives (see \ref{mot} and \ref{lef}). However, in the general framework of algebraic geometry we usually deal with singular varieties (or just open smooth schemes) and their cohomological invariants are associated with ``extensions'' (see \ref{ext}).

\subsubsection{ } \label{mm}  The general setting of  {\it mixed} motives takes care of such extensions classifying the cohomology theories defined on all algebraic varieties or schemes. Grothendieck \cite{GM} and Beilinson  (see Jannsen \cite{M} for a complete account\footnote{There are also some motivic conjectures on L-functions, see \cite{B}.}) conjectural theory  push further up to include this picture. Let $Sch_k$ now denote all algebraic varieties (or schemes) and imagine the corresponding cohomology theories (\cf \ref{ex}, \ref{mhs} and \ref{hdg}) carrying on mixed structures. The category $\cM_k^{m}$ of mixed motives over $k$ should then be an abelian $\Q$-linear category containing Grothendieck pure motives $\cM_k$. If $\cM_k^{m}$ exists we can consider the ``derived category'' $D^b(\cM_k^{m})$ of bounded complexes\footnote{Complexes and derived categories are a quite subtle technical device in order to treat homological algebra contructions. The main conceptual tasks are the ``cone'' of a map yielding a long exact sequence in homology and  the ``quasi-isomorphim'' inducing isomorphisms in homology.}  and we expect a motivic complex $M(X)\in D^b(\cM_k^{m})$ for $X\in Sch_k$ such that, for  $X\in \cV_k $ smooth and projective
$$M^i(X) = H^i(M(X))\in \cM_k$$ where $H^i$ is the homology of a complex. Moreover, the higher $\Ext^i$ in $\cM_k^{m}$ would have a geometrical meaning, see \ref{filt}. However, the current work on the construction of $\cM_k^{m}$ associated to arbitrary $k$-varieties and its full theory is mainly conjectural (see \cite{M} and \cite{BV} for the theory of  mixed motives of level $\leq 1$, \cf \ref{1-mot} below). A {\it motivic cohomology}\,  $\oplus_i M^i(X)[-i]\in \cM_k^{m}$, for arbitrary schemes $X\in Sch_k$, in Grothendieck sense should be at least universal among a well defined cohomology theory. In order to deal with cohomological motives of singular schemes we need a ``motivic'' description of the Picard functor  (see \cite{BV} and also Voevodsky's commentary \cite[p. 195]{V}, \cf \ref{der1m}). 

However, Voevodsky's homotopy invariant (homological) theory is suitable for several important purposes. Voevodsky's triangulated\footnote{Triangulated categories are an axiomatic version of derived categories.} category of geometrical motives over a field is also an attempt to construct a version of $D^b(\cM_k^{m})$ without assuming the existence of $\cM_k^{m}$.  We sketch the construction below.\footnote{We also have somewhat different constructions due to M. Levine and M. Hanamura. These approaches are based on the theory of algebraic cycles and the relations between them are quite well understood. Different versions of the triangulated category of motives over a field are historically the first examples of ``motivic homotopy categories''. The theory of these categories is closely related to the theory of homotopy invariant sheaves and cohomologies.} See \cite{V} and \cite[Part II]{A} for details.

\subsubsection{} \label{trim} Let $Sm_k$ denote the category of smooth algebraic schemes over $k$. For a pair $X,Y\in Sm_k$ we let $c (X, Y)\subseteq \cZ^*(X\times Y)$ denote the sub-group of {\it finite} correspondences: it is generated by sub-schemes $Z \subseteq X\times Y$ which are finite over $X$ and surjective on a connected component of $X$. A finite correspondence  $X\leadsto Y$  is somewhat a ``finite multivalued'' function from $X$ to $Y$. Similarly to \ref{cor} we get an additive category $SmCor (k)$ such that objects are smooth schemes of finite type over $k$ and morphisms are finite correspondences. Associating  a map $f:X\to Y$ to its graph $\Gamma_f$ we obtain
$[ - ]: Sm_k\to SmCor (k)$
where $[X]$ just denotes the object of $SmCor (k)$ corresponding to $X\in Sm_k$.
The triangulated category $\DM_{\gm}^{\eff}(k)$ of effective geometrical motives 
is the pseudo-abelian envelope\footnote{Formally adding kernels and cokernels of projectors as in \ref{mot} above.} of the localization of the homotopy category of bounded complexes  $K^b(SmCor (k))$ with respect to the thick sub-category generated by 
the following complexes:
\begin{enumerate}
\item {\it Homotopy}: for $X\in Sm_k$ we have $[X\times \Aff^1]\to [X]$ induced by the projection;
\item {\it Mayer-Vietoris}:  for $X= U\cup V\in Sm/k$ an open covering we get
$$[U\cap V]\to [U]\oplus [ V]\to [X]$$
\end{enumerate}
Denoting $M (X)$ the resulting object associated to $[X]$ we obtain a covariant functor $M: Sm_k \to \DM_{\gm}^{\eff}(k)$ from the category of
smooth schemes of finite type over $k$. We have then forced homotopy invariance by formally inverting (1) so that $$M (X\times \Aff^1) = M (X)$$ in $\DM_{\gm}^{\eff}(k)$. Moreover, from (2) we obtain a distinguished triangle
$$M (U\cap V)\to M(U)\oplus M (V) \to M (X)\longby{+1}$$
granting the Mayer-Vietoris axiom. There is a tensor structure such that $M (X\times Y)= M (X)\otimes M(Y)$.

\subsubsection{} \label{trimax}  Over a field $k$ of zero characteristic $M : Sch_k \to
\DM_{\gm}^{\eff}(k)$ extends to all schemes of finite type without change. On the other hand we have a functor $\cM^{\eff}_{\equiv_{\rm rat}}\to \DM_{\gm}^{\eff}(k)$ such that $M (X)$ is consistent with the description in \ref{cyc}-\ref{mot} when $X\in \cV_k$. 
In particular, the motive of a point $M (pt)\df \Z$ is the unit of the tensor structure and the structural morphism $X\to k$ yields the complex $[X]\to [k]$ and a triangle
$M_{\rm red}(X)\to M(X)\to \Z$. The Tate object is defined as $\Z (1)\df M_{\rm red}(\P^1)[-2]$  and $\DM_{\gm}(k)$ is defined by inverting $\Z (1)$.\footnote{Note that with respect to \ref{lef} here we are dealing with a homological theory. See \cite[II.17]{A} for a digression on this point.}

\subsubsection{} \label{motcompl}  Actually, Voevodsky also provides a larger category $\DM_{-}^{\eff}(k)$ of (effective) motivic complexes over $k$. This is suitable in order to reinterpret $M (X)$ {\it via}\, ``homotopy invariant sheaves with transfers''. Define the representable presheaf with transfers $L (X)$ on $SmCor (k)$, for $X\in Sm_k$, by $L (X) (Y)\df c(X, Y)$ yielding a sheaf  (for the Nisnevich topology) on $Sm_k$. A presheaf $\cF$ on  $SmCor (k)$ is homotopy invariant  if $\cF (X)=\cF (X\times \Aff^1)$ and the resulting category of sheaves $HI_k$ is an abelian category. The full sub-category  $\DM_{-}^{\eff}(k)\subset D^{-}(Shv_{Nis}(SmCor(k)))$ consists of complexes with homotopy invariant cohomology sheaves.\footnote{The canonical t-structure induces the so called homotopy t-structure on $\DM_{-}^{\eff}(k)$ with heart $HI_k$.} If  $\cF$ is a sheaf with transfers we can apply the Suslin complex $C_*$ given by $C_i(\cF)(Y) \df \cF (Y\times \Delta^i)$
and we obtain a complex of sheaves with transfers with homotopy invariant cohomology sheaves. For $\cF = L(X)$ the complex $C_*(L (X))\df C_*(X)$ is called the Suslin complex of $X$ and provides an algebraic version of singular homology.
One of the main technical achievements due to Voevodsky is the full embedding of  $\DM_{\gm}^{\eff}(k)$ into $\DM_{-}^{\eff}(k)$ in such a way that the image of $M(X)$ is $C_*(X)$. 

\subsubsection{} \label{motcom}  From the above discussion we can get a formula 
$$\Hom_{\DM_{\gm}^{\eff}} (M (Y),M(X)[j])= \HH^j(Y,C_*(X))$$
for $X,Y\in Sm_k$ and the following groups $H_j^s(X,\Z)\df \Hom_{\DM} (\Z[j],M(X))$ are the so called  Suslin homology groups. For example, Suslin-Voevodsky have shown that they agree with the singular homology of $X_{\rm an}$ when both are taken with finite coefficients.
Another useful task of  $\DM_{-}^{\eff}(k)$ is that it has ``internal'' Hom-objects $\ihom$ for morphisms from objects of $\DM_{\gm}^{\eff}$. Moreover, Voevodsky defined motivic complexes $\Z (j)$ and  corresponding ``motivic cohomology'' {\sf groups} by 
$$H^i_m (X,\Z (j))\df \Hom_{\DM} (M(X), \Z(j)[i])$$
Remark that this ``cohomology'' is not of the kind listed in \ref{ax} or \ref{ext}.
In particular, here $\Z (1)=\G_m[-1]$ as homotopy invariant sheaf with transfers and we obtain $H^2_m (X,\Z (1))=\Pic (X)= CH^1 (X)$ for $X$ smooth. Remark that if $X$ is not smooth $\Pic (X)$ is not homotopy invariant. Moreover, there are \'etale versions of motivic complexes  $\Z /n (j)$ whose hypercohomology coincide with $H^i(X_{\et}, \mu_{n}^{\otimes j})$ where $\mu_{n}$ is the \'etale sheaf of $n$-th roots of unity\footnote{For example, if $k =\bar k$ is algebraically closed we always  
have a non-canonical isomorphism $\mu_{n}^{\otimes j}\cong \Z/n$ by choosing a primitive root of unity.} and $n$ is prime to the characteristic of $k$ (\cf \ref{ex}.(2)).
These motivic cohomology groups are also related to algebraic $K$-theory and have been employed in the proof of 
\begin{itemize}
\item the {\em  Kato conjectures} claiming that $K_i^M(k)/n = H^i(k_{\et}, \mu_{n}^{\otimes i}) $ where $K_i^M(k)$ denotes Milnor $K$-theory of the field $k$ and $(n, \car (k) )= 1$.
\end{itemize}
We redirect the interested reader to forthcoming specific articles on the proof of Kato conjectures and to the existing good survey on the Milnor conjecture \cite{K}. 

We finally, briefly, mention (see \cite[II.22]{A} for more details) that there are mixed realization functors, \eg to the derived category of mixed Hodge structures
$$R:\DM_{\gm}\to D^+({\rm MHS})$$
where, for simplicity, we omit reference to the coefficients. These functors induce homomorphisms from ``motivic cohomology'' groups to similarly defined ``cohomologies'', the so called {\it absolute}\,  cohomologies, \eg absolute Hodge cohomology, \cf \ref{ajdel}. There is a conjectural picture (see \cite{B}) regarding values of $L$-functions of motives involving such homomorphisms which are called {\it regulators}\, in the current terminology.\footnote{A remark of M. Saito: ``Roughly speaking, there are two kinds of cohomology. One is Betti, de Rham, singular, $\ell$-adic, etc. The second contains Deligne cohomology, the absolute Hodge cohomology, the absolute \'etale (or continuous) cohomology, and the motivic
cohomology groups, \ie (higher) Chow groups.
The relation with motives is that the motive is universal among the
cohomology in the first sense, \ie the cohomology of the first kind
is obtained by applying a `forgetful' functor to the motives.
However, the cohomology in the second sense is obtained from the
motives by using something more complicated, \eg the group of
morphisms in an appropriated category.

Note that the motivic cohomology groups are not universal among the cohomology
in the first sense, but should be `universal' among the cohomology in
the second sense. (Here universal would mean simply that there is a
canonical morphism from the motivic cohomology to the other cohomology
of the second kind.)
I do not know whether we can get for example the Betti cohomology from
the motivic cohomology (as a vector space without additional structure).

Some other difference is that Tate twist is defined only for the
cohomology in the first sense, or rather, Tate twist defines an
`automorphism'. For the second kind, it is usually indexed by two
indices, and one of them corresponds to the Tate twist.

The cohomology of the first kind may be called Weil-type
cohomology because they satisfy the axioms of Weil cohomology
(here the Tate twist is usually omitted because it gives an
`automorphism'). The cohomology of the second kind may be called Deligne-type
cohomology. For the latter the axioms of Weil cohomology are not
satisfied.''}

\section{Improvements}
\subsection{Weil Conjectures} The most ``classical application'' of Grothendieck motives \ref{mot} is to the third of the Weil conjectures which we briefly mention here. For a more accurate, still introductory, explanation see Kleiman \cite{M}. For $X$ (smooth and projective) defined over a finite field with $q$ elements, we have the action of the Frobenius $\sigma : X \to X$ (which carries a point $x$ to the point $x^q$ whose coordinates are the $q^{th}$ powers of $x$). For each $n\geq 1$ let $a_n$ be the number of points of $X$ with coordinates in a finite extension field with $q^n$ elements. This $a_n$ also equals the number of fixed points of the iterated Frobenius $\sigma^n$ on $X$ and we have that
$$Z (t) \df \exp \Bigg( \sum_{n\geq 1} \frac{a_n t^n}{n}\Bigg)$$
is the famous zeta function of $X$. The $\sigma^n$ induces an action on the cohomology groups $H^i (X)$. Moreover, for each $n\geq 1$ there is a {\it Lefschetz trace} formula
$$a_n = \sum_{i=0}^{2\dim (X)}(-1)^i{\rm Tr} (\sigma^n\mid_{H^i (X)})$$
where Tr denotes the trace of the induced action of $\sigma^n$ on each $H^i (X)$ (we here tacitly  deal with \'etale $\ell$-adic cohomology in \ref{ex}.(2)).
Therefore, we obtain that 
$$a_n = 1 + \sum_{i=1}^{2\dim (X)-1}(-1)^i \sum_j\lambda_{ij}^n+ q^{n\dim (X)}$$
where $\lambda_{ij}$ are the eigenvalues of $\sigma$ on $H^i (X)$. Now the third of
Weil's conjectures\footnote{The first is that $Z (t)$ is a rational function, which is a formal consequence of the Lefschetz trace formula, and the second is that there is a functional equation, which follows from Poincar\'e duality, so that Grothendieck $\ell$-adic cohomology theory suffices.} is that the characteristic polynomial of  $\sigma\mid_{H^i (X)}$ has integer coefficients, which are independent of the cohomology theory and the eigenvalues $\lambda_{ij}$ are of absolute value $q^{i/2}$. This is a quite trivial consequence of the Standard Conjectures.\footnote{The third of the Weil conjectures was actually proven by Deligne (1973) by a different method.} In fact, in general, let $\gamma: X\leadsto X$  be {\it any} correspondence defined by an algebraic cycle, then the characteristic polynomial of $\gamma : H^i (X)\to H^i (X)$ has integer coefficients which are independent of the cohomology theory, if the Standard Conjectures hold, since we have another {\it trace} formula
$${\rm Tr} (\gamma\mid_{H^i (X)}) = (-1)^i \gamma \cdot p_{2\dim (X)-i}$$
where $ p_{2\dim (X)-i}$ is the transposed of $p_i$ in \ref{gal} and $ \gamma\cdot p_{2\dim (X)-i}$ are the resulting intersection numbers. Moreover, the so called ``Betti numbers'',\ie  $\dim H^i(X)$, are independent of the cohomology theory (to see this one takes $\gamma =1$ in the trace formula).

\subsection{Deligne Conjectures} 
A classical construction in algebraic geometry associates to a smooth projective curve $X$ its Jacobian variety $J(X)$. This is a so called ``abelian variety'' which is a smooth projective algebraic variety along with a structure of abelian group on its points. 
For such a curve $X$ there is a (pointed) morphism to the Jacobian $X\to J (X)$ such that any (pointed) morphism from $X$ to an abelian variety factors through $J(X)$, \ie $ J (X)$ can  also be characterized by this universal property.  This is a key point in algebraic geometry and can be extended to higher dimensional smooth projective $X$ so that the resulting universal abelian variety $\Alb (X)$ is the so called {\it Albanese variety} and the map $X\to \Alb (X)$ the Albanese map.

The $\Q$-linear category of abelian varieties up to isogeny\footnote{An isogeny between abelian varieties is a surjective morphism with finite kernel.} is an abelian semi-simple category  (which can be obtained as the pseudo-abelian envelope of the category of Jacobians and $\Q$-linear maps). Recall \ref{lef} that in the case of curves $M^1 (X)\in \cM^{\eff}_{\equiv_{\rm rat}}$ is the motive of $X$ (with $\Q$-coefficients)  refined from lower and higher trivial components, such that, for smooth projective curves $X$ and $Y$ we have the following nice formula
$$\Hom (M^1(X),M^1(Y)) \cong \Hom (J(X), J (Y))_{\sQ}$$ due essentially to Weil  (see \cite{MU}). Thus, as pointed out by Grothendieck, the theory of {\em pure} motives of smooth projective curves is equivalent to the theory of abelian varieties up to isogeny.

\subsubsection{}\label{1-mot} Starting from the generalization of the Jacobian of a curve
 Deligne \cite{D} provided a theory of  (free) {\it 1-motives} and, in such a way, also provided the motivic cohomology of  possibly singular curves. This theory is workable by making algebraic the definition of mixed Hodge structures of level $\leq 1$ and then yields a corresponding theory of mixed motives of level $\leq 1$. A {\it free} 1-motive over a field $k$ is a two terms complex  $L\to G$ where $L$ is a finitely generated {\it free} abelian group and $G$ is an extension of an abelian variety by a torus.  The category $\cM_1^{\fr}(k)$ of 1-motives over $k$ has objects $M\df [L\to G]$ and morphisms are pairs of maps making a commutative square. For $k=\C$ we have an equivalence of categories
 $$\cM_1^{\fr}(\C)\cong {\rm MHS}_1^{\fr}$$
 with the full sub-category  ${\rm MHS}_1^{\fr}\subset  {\rm MHS}$ of torsion-free (graded polarizable)  mixed Hodge structures of level $\leq 1$. Remark that Deligne's category $\cM_1^{\fr}(k)$ is not abelian. In order, to get an abelian category we have to allow {\it torsion} 1-motives taking as objects $M\df [L\to G]$ where now $L$ may have torsion and formally invert quasi-isomorphisms $M\to M'$. The resulting category 
 $\cM_1(k)$ is abelian and 
$$\M (\C)\cong {\rm MHS}_1$$ is an equivalence with (graded polarizable)  mixed Hodge structures of level $\leq 1$. See \cite{BV} for a full account and references.

\subsubsection{}\label{dconj} For higher dimensional varieties Deligne \cite{D} proposed some conjectures which imply the algebraic nature of certain complex tori obtained by trascendental methods as follows.

Let $X$ be a complex algebraic variety of dimension $\leq n$ and let $H^*(X, \Z )$ be the mixed Hodge structure on the singular cohomology of the associated analytic space.
Denote $H^*_{(1)}(X, \Z )\subseteq H^*(X, \Z )$ the largest sub-structure and $H^*(X, \Z)^{(1)}$ the largest quotient of level $\leq 1$. Now, further deleting torsion, since ${\rm MHS}_1^{\fr}\cong\cM_1^{\fr}(\C)$ we obtain corresponding 1-motives over $\C$  and 
\begin{itemize}
\item  the {\it Deligne conjectures} claim that these 1-motives are algebraically defined over any field $k$, \ie when $k=\C$, the mixed Hodge structures
$H^i_{(1)}(X, \Z (1))_{\fr}$, $H^i(X, \Z (i))^{(1)}_{\fr}$ for $i\leq n$ and
$H^i(X, \Z (n))^{(1)}_{\fr}$ for $i\geq n$ are provided by algebraic  methods only.\footnote{Here $ \Z (\dag)$-coefficients denote that the mixed Hodge structure is Tate twisted by $\dag$.}
\end{itemize}
For example (see \cite{BV} for references and more details on these conjectures) we can construct 1-motives with torsion $\Pic^+(X,i)\in \M(k)$ proving the conjecture for $H^i_{(1)}(X, \Z (1))_{\fr}$ up to isogeny, \eg for $X$ a smooth projective curve $\Pic^+(X,0)= J (X)$ corresponds to $M^1 (X)$ and, over $\C$, to $H^1(X, \Z )$.
We also have ``Cartier duals'' (up to isogeny) $\Alb^-(X,i)$ such that $\Alb^-(X,0)=\Alb (X)$ is the classical Albanese variety for $X$ smooth.

\subsubsection{}\label{der1m} Since a 1-motive $M =[L\to G]$ is actually a complex of \'etale sheaves, where $L$ and $G$ are clearly homotopy invariants, it yields an effective complex of homotopy invariant \'etale sheaves with transfers, and hence an object of $\DM_{-,\et}^{\eff}(k)$ (\cf \ref{motcompl} and \cite{BV}). Considering the $\Q$-linear category $\M^{\sQ}$ of 1-motives up to isogeny we obtain a functor $\M^{\sQ}\to\DM_{-}^{\eff}(k;\Q)$ and
we can see that there is a fully faithful functor $$\Tot : D^b(\M^{\sQ}) \by{\simeq} d_{\leq 1}\DM_{\gm}^{\eff}(k)\subseteq \DM_{\gm}^{\eff}(k;\Q)$$ whose essential image is  the thick triangulated sub-category $d_{\leq 1}\DM_{\gm}^{\eff}(k)\subseteq \DM_-^{\eff}(k)$ generated by motives of smooth varieties of dimension $\leq 1$ (see \cite{V} and \cite{A}). We can show (in a joint work with B. Kahn) that there is a functor
$$\LAlb : \DM_{\gm}^{\eff}(k;\Q)\to D^b(\M^{\sQ})$$
which is left adjoint to the embedding $\Tot $ above (see \cite{BV}). 
The aim is that this operation will be rendering the 1-motives predicted by the Deligne conjecture. When applied to $M(X)$, denote $\LAlb(X)\df \LAlb(M (X) )$ and thus obtain a universal map $M (X)\to \Tot \LAlb(X)$ in $\DM^{\eff}_{\gm}(k;\Q)$ (which is an isomorphism if $\dim (X)\leq 1$). This motivic Albanese map is refining the classical Albanese map. 
The motivic Albanese complex $\LAlb(X)$ is a bounded complex of 1-motives and the 1-motivic homology for $i\geq 0$ $$\LA{i}(X)\df H^{-i}(\LAlb(X))$$ will coincide with
$\Alb^-(X, i-1)$ (for $i\geq 1$ up to isogeny). 

Dually, denote $\RPic (X)\df \LAlb(X)^{\vee}$ and $\RA{i} (X) = H^i (\RPic (X))$. 
This $\RPic (X)$ is the motivic Picard complex. For example, we expect that (the homotopy invariant part of) the first motivic cohomology $M^1 (X)\in \cM^m_k$ (if it exists!) is given by $\RA{1} (X) \in \cM_1\subset  \cM^m_k$ for an arbitrary $X\in Sch_k$ (\cf \ref{mm} and \ref{2-mot}).

\subsection{Bloch Conjectures} 
Recall that the Albanese map $X\to \Alb (X)$ is also yielding a (surjective) map $CH^n(X)_{\deg 0}\to \Alb (X)$ from the Chow group of points of degree zero on $X$ smooth projective $n$-dimensional. Based on an argument due to Severi, Mumford was able to show that the kernel is huge if there is a non trivial holomorphic $2$-form on a surface. On the contrary we have
\begin{itemize}
\item  the {\it Bloch conjecture} that, for a surface $X$, the Albanese kernel $\ker (CH^2(X)_{\deg 0}\to \Alb (X))$ vanishes if $H^2(X_{\rm an},\Z)$ is algebraic (\ie $A^1 =H^2$ or equivalently if  all global holomorphic $2$-forms vanish).
\end{itemize}
The named condition on 2-forms is also equivalent to the equality $N^1H^2(X)= H^2(X)$ given by the coniveau filtration \ref{ghdg}.

\subsubsection{}\label{ajdel} Here we have a general picture for $X$ a compact complex algebraic manifold by making use of the Deligne cohomology
$H^*(X,\Z(\cdot)_{\cD})$ which is defined by taking the hypercohomology of $\Z \to \Omega^{<\cdot}_X$, \ie the truncated De Rham complex augmented over $\Z$ (see \cite{B}). There is a forgetful map $H^*(X,\Z(\cdot)_{\cD})\to H^*(X_{\rm an},\Z(\cdot))$ and there is a cycle map $$c\ell_{\cD} : CH^i(X)\to H^{2i}(X,\Z(i)_{\cD})$$ yielding by composition the usual cycle map \ref{cl}. Moreover, we have an extension 
$$0\to J^i(X)\to H^{2i}(X,\Z(i)_{\cD})\to H^{i,i}_{\sZ}(X)\to 0$$
where $J^i(X)$ is the Griffiths intermediate Jacobian\footnote{This is a complex torus  given by the $\C$-vector space $H^{2i-1}_{DR}(X)/F^i$ modulo the image of $H^{2i-1}(X,\Z)$.} and the cycle map restricts to the so called Abel-Jacobi map $CH^i(X)_{\rm hom}\to J^i(X)$ where $CH^i(X)_{\rm hom}\df \ker ( CH^i(X)\to H^{2i}(X_{\rm an},\Z(i))$. In particular, $CH^n(X)_{\deg 0}= CH^n(X)_{\rm hom}$ and $J^n(X) = \Alb (X)$ for $X$ $n$-dimensional so that the ``Albanese kernel'' equals the kernel of $c\ell_{\cD}$.
Moreover,  for a surface $X$ as above, such that $H^2(X,\Z)$ is algebraic,
we can show that the ``Albanese kernel'' vanishes if and only if for any
non-empty Zariski open $U\subset X$ and $\omega\in
H^2(U,\C)$, there is a non-empty Zariski open
$V\subset U$ such that $$\omega\mid_V =
\omega'+\omega_{\sZ},$$ where $\omega'\in F^2H^2(V,\C)$, \ie a meromorphic $2$-form with logarithmic poles, and $\omega_{\sZ}$ is integral, \ie $\omega_{\sZ}\in\im H^2(V,\Z (2))$. This is a quite explicit but weird unproven property of local complex 2-cohomology classes. 

\subsubsection{}\label{filt}  Regarding the $\text{``Albanese kernel''} \subseteq CH^n(X)_{\rm hom}\subseteq CH^n(X)$ Bloch \cite{BL} observed that this is somewhat a filtration on $CH^n(X)$ of a $n$-dimensional $X$ which should be detected by motivic arguments, \ie by the action of correspondences. In general, Murre \cite{MU} conjectured a motivic filtration $F^*_m$ on  $CH^i(X)$ given by a lifting $\pi_i$ of the K\"unneth components $p_i$ of the diagonal cycle, \ie  such that $F^1_m\df \ker \pi_{2i}\mid_{CH^i}$, $F^2_m\df \ker \pi_{2i-1}\mid_{F^1}$, \ldots ,$F^{i+1}_mCH^i =0$. Assuming the existence of mixed motives $\cM^m_k$ it can be resumed by Beilinson conjectural formula (see Jannsen \cite{M}) 
$$ \gr_{F}^j CH^i(X) = \Ext_{\cM^m_k}^j(\Q, M^{2i-j}(X)(i))$$
formally given by the Ext spectral sequence\footnote{The spectral sequence is $E^{p,q}_2=\Ext^p(\Q , M^q(X)(\cdot))\implies \EExt^{p+q}(\Q , M(X)(\cdot))$.} associated to the bounded motivic complex $M (X)\in D^b(\cM^m_k)$.
Here $F^0_mCH^i(X) = CH^i(X)$, $F^1_mCH^i(X) = CH^i(X)_{\rm hom}$ and
$$ \gr_{F}^0CH^i(X) = \Ext_{\cM^m_k}^0(\Q, M^{2i}(X)(i))= \Hom_{\cM_k}(\Q , M^{2i}(X)(i))$$
holds by construction of Grothendieck motives $\cM_k\subseteq \cM^m_k$. Moreover one obtains that
 $F^2_mCH^n(X) = \text{``Albanese kernel''}$ on $X$ $n$-dimensional and 
 $$ \gr_{F}^1CH^n(X) = \Ext_{\cM^m_k}^1(\Q, M^{2n-1}(X)(n))= \Alb (X)$$
For $n=2$, \ie on a surface $X$ we have $F^{3}_mCH^2 =0$, we thus obtain
 $$F^2_mCH^2(X) = \text{``Albanese kernel''}=  \Ext_{\cM^m_k}^2(\Q, M^{2}(X)(2))$$ where $M^{2}(X)= (X, p_2)\in \cM_k$, therefore (exercise!): the Albanese kernel vanishes if $H^2(X)$ is algebraic.

\subsubsection{}\label{2-mot} Finally remark that a main step here will just be the construction of 2-motives $\mathcal{M}_2$ containing the 1-motives $\mathcal{M}_1$ in \ref{1-mot} and, by the way, further up to $n$-motives $\mathcal{M}_n(k)\subseteq \cM^m_k$ as targets for motives of algebraic $k$-schemes of dimension $\leq n$. A decomposition $M(X) =\oplus M^i(X)[-i]$ of the  motivic complex of an algebraic variety $X$ will be refined by asking $M^i(X)\in \mathcal{M}_i$ providing the motivic cohomology of such $k$-schemes. Such $\mathcal{M}_i(k)$ will be abelian categories such that $\Ext^*_{\mathcal{M}_i}=0$ for $*>i$ when $\Q$-linearized or $k=\bar k$. There will be a compatible family of Hodge realizations as in the following conjectural picture
$$\begin{array}{ccc}
\mathcal{M}_i(\C)&\longby{R_i}& \mbox{MHS}_i\\
\downarrow & &\downarrow\\
\cM^m(\C)&\longby{R}& \mbox{MHS}
\end{array}$$
Here ${\rm MHS}_i$ will be an abelian category generalizing mixed Hodge structures along with a forgetful functor ${\rm MHS}_i\to {\rm MHS}$ with essential image mixed Hodge structures of level $\leq i$.  The wish is now that the Hodge realization $R$ of the motive of an algebraic variety which is of level $\leq i$ would be lifted back {\it via}\, $R_i$ and algebraically defined in $\mathcal{M}_i$. For example, if $H^2(X_{\rm an},\Z)$ is algebraic for a surface $X$ then all $H^*(X_{\rm an},\Z (1))\in \mbox{MHS}_1$ is of level $\leq 1$ and Bloch's conjecture is a consequence of the following property\\

\begin{itemize}
\item  if $H^*(X_{\rm an},\Z (1))\in \mbox{MHS}_1$ then $M (X)\in \mathcal{M}_1$ is 1-motivic.\\
\end{itemize}

\end{document}